\documentclass[12pt]{article}

\usepackage{amsmath}
\usepackage{mathrsfs}
\usepackage{bbm}
\usepackage{amssymb,a4}
\usepackage{amsfonts,amssymb,amsthm}
\usepackage{extarrows}
\usepackage{color}
\usepackage[all,pdf]{xy}
\usepackage{multicol}
\usepackage{fancyhdr}

\allowdisplaybreaks  % 允许过长的公式分页

\newtheorem{thm}{\bf Theorem}[section]

\newtheorem{lem}[thm]{\bf Lemma}

\newcommand{\Z}{{\mathbb Z}}
\newcommand{\C}{{\mathbb C}}
\newcommand{\N}{{\mathbb N}}
\newcommand{\zp}{{\Z_+}}
\newcommand{\g}{\mathfrak{g}}
\newcommand{\U}{{\mathcal U}}
\newcommand{\al}{\alpha}
\newcommand{\e}{{\epsilon}}
\newcommand{\p}{\partial}
\newcommand{\sgm}{\sigma}

\newcommand{\ov}{\overline}
\newcommand{\spanc}[1]{\mathrm{span}_\C\left\{#1\right\}}

\newcommand{\nd}{\N^d}
\newcommand{\ndz}{\nd\backslash\{\bff 0\}}
\newcommand{\ad}{\mathrm{ad}}
\newcommand{\zd}{\Z^d}
\newcommand{\cd}{\C^d}
\newcommand{\cq}{\C_Q}

\newcommand{\D}{\mathcal{D}}
\newcommand{\calg}{\mathcal{G}}
\newcommand{\calgx}{\calg^{(x)}}
\newcommand{\calgt}{\calg^{(t)}}
\newcommand{\vavw}{{\mathcal V^\al(V,W)}}
\newcommand{\ZZ}{\mathcal{Z}}
\newcommand{\ZD}{\ZZ\D}

\newcommand{\mk}[1]{M_{k_{#1}}(\C)}
\newcommand{\gl}{\mathfrak{gl}}
\newcommand{\gln}{\mathfrak{gl}_N}
\newcommand{\gld}{\mathfrak{gl}_d}
\newcommand{\bff}[1]{{\bf #1}}
\newcommand{\ltwo}[2]{L({{\bf #1},{\bf #2}})}  %L(m,r)
\newcommand{\sgmf}[2]{\sgm({\bf #1},{\bf #2})}% sigma function
\newcommand{\parder}[1]{\frac{\p}{\p {#1}}}     %\partial derivatives
\newcommand{\partwo}[2]{\p(\bff #1,\bff #2)}     %\partial(u,m)
\newcommand{\parthree}[3]{\p(\bff #1,\bff #2+\bff #3)}

\newcommand{\tone}[1]{t^{\bf #1}}    % monomomial in t
\newcommand{\ttwo}[2]{t^{{\bf #1}+{\bf #2}}}
\newcommand{\tthree}[3]{t^{{\bf #1}+{\bf #2}+{\bf #3}}}
\newcommand{\tfour}[4]{t^{{\bf #1}+{\bf #2}+{\bf #3}+{\bf #4}}}
\newcommand{\ovbff}[1]{\ov{\bff #1}}
\newcommand{\tonebar}[1]{t^{\bf\ov{#1}}}

\newcommand{\xone}[1]{x^{\bf #1}}    % monomomial in x
\newcommand{\xtwo}[2]{x^{{\bf #1}+{\bf #2}}}
\newcommand{\bxone}[1]{X^{\bf #1}}    % monomomial in big X
\newcommand{\bxtwo}[2]{X^{{\bf #1}+{\bf #2}}}
\newcommand{\ftwo}[2]{f_{\bff #1}(\bff #2)}
\newcommand{\gtwo}[2]{g_{\bff #1}(\bff #2)}

\newcommand{\done}[1]{d_{\bff #1}}
\newcommand{\dtwo}[2]{D({\bff #1},{\bff #2})}

\newcommand{\inr}{\in R}
\newcommand{\notinr}{\notin R}

\newcommand{\pf}[1]{\noindent{\bf Proof.}#1\hfill{}$\Box$}

\title{\bf\Large Cuspidal modules for the derivation Lie algebra\\
 over a rational quantum torus}
\author{Chengkang Xu$^{\ast}$\\
{\small School of Mathematical Science,
Shangrao Normal University, Shangrao, China
}}

\date{}

\begin{document}
\maketitle
\renewcommand{\thefootnote}

\setcounter{footnote}{-1}\footnote{* Corresponding author.
The author is supported by the National Natural Science Foundation of
China(No. 11626157, 11801375).
\emph{E-mail:} xiaoxiongxu@126.com.}

{\noindent\bf \normalsize Abstract:}  {\small
Let $\cq$ denote a rational quantum torus with $d$ variables,
and $\ZZ$ be the centre of $\cq$.
In this paper we give a explicit description of the structure of the cuspidal modules
for the derivation Lie algebra $\D$ over $\cq$,
with an extra associative $\ZZ$-action.
%As a consequence, we prove that any irreducible cuspidal $\dercq$-module
%with an extra associative $\ZZ$-action is isomorphic to $\vavw$,
%for some $\al\in\cd$, finite dimensional irreducible $\gld$-module $V$
%and finite dimensional $\Gamma$-graded $\gln$-module $W$.
%This was originally proved in \cite{LiuZ2}, but in a more computational way.
%Moreover, from our proof one can see clearly
%how the mysterious algebra $\gln$ and the $\Gamma$-graded $\gln$-module come into the picture,
%and their relation with the $\gld$-module $V$.
}

\medskip
{\noindent \small {\bf Keywords}: }{\small derivation Lie algebra; cuspidal module;
rational quantum torus; module of tensor fields.}

%%%%%%%%%%%%%%%%%%%%%%%%%%%%%%%%%%%%%%%%%%%%%%%%%%%%%%%%%%%%%%%%%%%%%%%%%%%%%%%%% Section 1
\section{Introduction}
\def\theequation{1.\arabic{equation}}
\setcounter{equation}{0}

In the past few decades, the representation theory of the derivation algebra
over a torus has attracted many mathematicians and physicists.
Let $d>1$ be a positive integer.
Denote by $A=\C[x_1^{\pm1},\cdots, x_d^{\pm1}]$
the commuting torus in variables $x_1,\cdots, x_d$,
and by $W_d$ the derivation algebra over $A$.
Irreducible modules with finite dimensional weight spaces over $W_d$
were classified in \cite{BF1}.
They are modules of tensor fields constructed independently
by Shen\cite{S} and Larsson\cite{Lar},
and modules of highest weight type given in \cite{BB}.
For the derivation algebra $\D$ over a rational quantum torus $\cq$, although
the modules of tensor fields were constructed in \cite{LT} and \cite{LiuZ2}
and the modules of highest weight type were given in \cite{Xu},
we are still far from the complete classification of irreducible modules
with finite dimensional weight spaces for $\D$.
However, some partial classifications are known to us.
For example, irreducible $\D$-modules with certain $\cq$-action were classified in \cite{RBS}.
In \cite{LiuZ2}, Liu and Zhao classified irreducible cuspidal $\D$-modules
with an extra associative action of the centre $\ZZ$ of $\cq$.
These modules were constructed using
a finite dimensional irreducible $\gld$-module $V$
and a finite dimensional $\zd/R$-graded irreducible $\gln$-module $W$.
Here $N$ is a positive number closely related to the structure of $\cq$,
and $R$ is a subgroup of $\zd$ corresponding to the center $\ZZ$.
To achieve their classification result in \cite{LiuZ2},
Liu and Zhao used heavy computations(totally five pages of them).
Moreover, the appearance of the algebra $\gln$ in the proof seems a little farfetched,
and it is not clear what is the relation between the $\gld$-module $V$
and the $\zd/R$-graded $\gln$-module $W$.

In this paper we study $\D$-modules(not necessarily irreducible)
with an extra associative $\ZZ$-action.
We will simply call such a module a $\ZD$-module.
Our main result is that
cuspidal $\ZD$-modules with support lying in one coset $\al+\zd, \al\in\cd$,
are in a one-to-one correspondence with finite dimensional $\Gamma$-graded
modules over a subalgebra $\calg$ of $\mathrm{Der}(AC)$,
where $AC=\C[x_1,\cdots,x_d]\oplus \cq/\ZZ$.
The algebra $\calg$ has a quotient isomorphic to $\gld\oplus\gln$.
We will use some results about the solenoidal Lie algebras over $A$ and $\cq$,
which were given in \cite{BF2} and \cite{Xu2} separately.
As a corollary, we reprove the classification of irreducible cuspidal $\ZD$-modules,
in a more conceptional way.
From our proof one can see clearly that the algebras $\gld, \gln$
and their modules appear naturally.
We mention that Liu and Zhao classified irreducible cuspidal $\D$-modules
with no other restriction
in \cite{LiuZ} using the classification result in \cite{LiuZ2}
and a method introduced in \cite{BF1}.

The present paper is arranged as follows.
In Section 2 we recall some results about the quantum torus $\cq$, the algebras $\D$,
the solenoidal Lie algebra $W_\mu$ over $A$
and the solenoidal Lie algebra $\g_\mu$ over $\cq$.
Section 3 is devoted to finite dimensional $\calg$-modules.
In Section 4 we prove our main theorem about the cuspidal $\ZD$-modules
and classify irreducible cuspidal $\ZD$-modules once again.

Throughout this paper, $\C,\Z,\N,\zp$ refer to the set of complex numbers,
integers, nonnegative integers and positive integers respectively.
For a Lie algebra $\mathcal{G}$, we denote by $\U(\mathcal G)$
the universal enveloping algebra of $\mathcal{G}$.
Fix $1<d\in \zp$ and a standard basis $\e_1,\cdots,\e_d$ for the space $\cd$.
Denote by $(\cdot\mid \cdot)$ the inner product on $\cd$.

%%%%%%%%%%%%%%%%%%%%%%%%%%%%%%%%%%%%%%%%%%%%%%%%%%%%%%%%%%%%%%%%%%%%%%%%%%%%%%%%%%%%% Section 2
\section{Notations and Preliminaries}
\def\theequation{2.\arabic{equation}}
\setcounter{equation}{0}

For a Lie algebra,
a weight module is called a {\bf cuspidal} module if all weight spaces are uniformly bounded.
The {\bf support} of a weight module is defined to be the set of all weights.

For $\bff m=(m_1,\cdots,m_d)^T\in\zd$ we denote by $\xone m$
the monomial $x_1^{m_1}\cdots x_d^{m_d}$ in $A$.
Then the derivation Lie algebra $W_d$ over $A$ has a basis
$\{\xone mx_i\parder x_i\mid \bff m\in\zd, 1\leq i\leq d\}$
subject to the Lie bracket
$$[\xone mx_i\parder x_i,\xone nx_j\parder x_j]=\xtwo mn (n_ix_j\parder x_j-m_jx_i\parder x_i).$$
Let $\mu=(\mu_1,\cdots,\mu_d)^T\in\cd$ be generic,
which means that $\mu_1,\cdots,\mu_d$ are linearly independent
over the field of rational numbers.
The solenoidal Lie algebra $W_\mu$ over $A$ is
the subalgebra of $W_d$ spanned by
$\{\xone m\sum_{i=1}^d\mu_ix_i\parder x_i\mid \bff m\in\zd\}$.
Many others call $W_\mu$ (centerless) higher rank Virasoro algebra.

Let $Q=(q_{ij})$ be a $d\times d$ complex matrix
with all $q_{ij}$ being roots of unity and satisfying
$q_{ii}=1,\ q_{ij}q_{ji}=1\text{ for all }1\leq i\neq j\leq d.$
The rational quantum torus relative to $Q$ is the unital associative algebra
$\cq=\C[t_1^{\pm1},\cdots, t_d^{\pm1}]$ with multiplication
$$t_it_j=q_{ij}t_jt_i\text{ for all }1\leq i\neq j\leq d.$$
For ${\bff m}=(m_1,\cdots, m_d)^T\in\zd$
we denote $\tone m=t_1^{m_1}\cdots t_d^{m_d}$.
For $\bff m, \bff n \in \zd$, set
$$\sgmf mn=\prod_{1\leq i<j\leq d}q_{ji}^{n_jm_i}\text{ and }
R=\{\bff m\in\zd\mid \sgmf mn=\sgmf nm\text{ for any }\bff n\in\zd\}.$$
Then the center $\ZZ$ of $\cq$ is spanned by $\{\tone m\mid\bff m\in R\}$.
By \cite{BGK}, the derivation Lie algebra $\D$ over $\cq$ has a basis
$$\{\tone m\p_i,\ad\tone s\mid\bff m\inr,\bff s\notinr\},$$
where $\p_i$ is the degree derivation with respect to $t_i$
defined by $\p_i\tone r=r_i\tone r$,
and $\ad\tone s$ is the inner derivation given by $\ad\tone s\cdot \tone r=[\tone s,\tone r]$.
For convenience we will simply write $\tone s=\ad\tone s$.

By Theorem 4.5 in \cite{N}, up to an isomorphism of $\cq$,
we may assume that
$q_{2i,2i-1}=q_i, q_{2i-,2i}=q_i^{-1}$ for $1\leq i\leq z$,
and other entries of $Q$ are all 1,
where $z\in\zp$ with $2z\leq d$
and the orders $k_i$ of $q_i$ as roots of unity satisfy
$k_{i+1}\mid k_i$ for $1\leq i\leq z$.
In this paper we always assume that $Q$ is of this form,
and keep the notation $q_i,z$ and $k_i$.
Then the subgroup $R$ of $\zd$ has the simple form
$$R=\bigoplus_{i=1}^z\left(\Z k_i\e_{2i-1}\oplus \Z k_i\e_{2i}\right)
\oplus \bigoplus_{l>2z}\Z\e_l,$$
and the center $\ZZ$ becomes a commuting torus on variables
$$t_1^{k_1},t_2^{k_1},t_3^{k_2},t_4^{k_2},\cdots,t_{2z-1}^{k_z},t_{2z}^{k_z},
t_{2z+1},\cdots, t_d.$$
Moreover, we have
$$\sgmf mr=\sgmf rm,\ \ \tone m\tone r=\ttwo mr\text{ for all }\bff m\in R, \bff r\in\zd.$$
Set $\partwo um=\sum_{i=1}^du_i\tone m\p_i$
for $\bff u=(u_1,\cdots,u_d)^T\in\cd$ and $\bff m\inr$.
Then the Lie bracket of $\D$ is
$$\begin{aligned}
 &[\partwo um,\partwo vn]=(\bff u\mid\bff n)\parthree vmn-(\bff v\mid\bff m)\parthree umn; \\
 &[\partwo um,\tone s]=(\bff u\mid\bff s)\ttwo ms; \ \ \
  [\tone r,\tone s]=(\sgmf rs-\sgmf sr)\ttwo rs,
 \end{aligned}$$
where $\bff u,\bff v\in\cd,\bff m,\bff n\inr$ and $\bff r,\bff s\notinr$.

Denote by $\mathcal I$ the ideal of the associative algebra $\cq$
generated by elements $\ttwo nr-\tone r, \bff n\in R,\bff r\in\zd$.
By \cite{N}, we have
$\cq/\mathcal I\cong \bigotimes_{i=1}^z\mk i\cong M_N(\C)$,
where $N=\prod_{i=1}^zk_i$ and
$M_n(\C)$ denotes the associative algebra of all $n\times n$ complex matrices.
It is well known that $\mk i$ can be generated as an associative algebra by
$$\begin{aligned}
&X_{2i-1}=E_{1,1}+q_iE_{2,2}+\cdots+q_i^{k_1-1}E_{k_i,k_i},\\
&X_{2i}=E_{1,2}+E_{2,3}+\cdots+ E_{k_i-1,k_i}+E_{k_i,1},
\end{aligned}$$
where $E_{kl}$ is the $k_i\times k_i$ matrix with 1 at the $(k,l)$-entry
and 0 elsewhere.
Denote $X^{\bff n}=\bigotimes_{i=1}^z X_{2i-1}^{n_{2i-1}}X_{2i}^{n_{2i}}$
for $\bff n\in\zd$.
It is easy to see that for any $\bff n\in R$, $\bxone n$ is the identity matrix in $M_N(\C)$,
and $\bxone r\bxone s=\sgmf rs\bxtwo rs$ for any $\bff r,\bff s\in\zd$.
Then the general linear Lie algebra $\gln=M_N(\C)$ has Lie bracket
$$[\bxone m,\bxone n]=(\sgmf mn-\sgmf nm)\bxtwo mn.$$
Set $\Gamma=\zd/R$ and let $\ovbff n$ denote the image of $\bff n\in\zd$.
Clearly, $\Gamma$ has order $N$
and the algebra $\gln$ has a $\Gamma$-gradation
$\gln=\bigoplus_{\ovbff n\in\Gamma}(\gln)_{\ovbff n}$
where $(\gln)_{\ovbff n}=\C\bxone n$.
In this paper by a $\Gamma$-gradation on $\gln$ we always mean this gradation.
Set
$$\Gamma_0=\{\bff n\in\zd\mid 0< n_{2i-1}\leq k_i,0< n_{2i}\leq k_i,1\leq i\leq z,
\text{ and }n_l=0,2z<l\leq d\},$$
which is a complete set of representatives for $\Gamma$.
When a representative of $\ovbff n$ is needed we always choose
the one $\bff n\in\Gamma_0$.

Let $\mu\in\cd$ be generic.
A solenoidal Lie algebra $\g_\mu$ over $\cq$ is the subalgebra of $\D$ spanned by
$$\{\partwo \mu m,\ \tone s\mid\bff m\inr,\bff s\notinr \}.$$
We mention that the subalgebra $\D_R$ spanned by
$$\{\partwo um\mid \bff u\in\cd,\bff m\inr\}$$
is isomorphic to the Lie algebra $W_d$ through the map defined by
$$ \partwo um\mapsto \xone n\sum_{i=1}^du_ik_ix_i\parder{x_i},$$
where $\bff m=B\bff n$ and $B$
is the $d\times d$ matrix $\mathrm{diag}\{k_1,k_1,k_2,k_2,\cdots,k_z,k_z,1,\cdots,1\}$.

%%%%%%%%%%%%%%%%%%%%%%%%%%%%%%%%%%%%%%%%%%%%%%%%%%%%%%%%%%%%%%%%%%%%%%%%%%%%%%%%%%%%% Section 3
\section{Finite dimensional $\calg$-modules}
\def\theequation{3.\arabic{equation}}
\setcounter{equation}{0}

In this section we introduce the Lie algebra $\calg$ and study its modules of finite dimension.
For $\bff s\notinr$ we will denote by $\tonebar s$ the image of $\tone s$ in $\cq/\ZZ$.
Let $AC=\C[x_1\cdots,x_d]\oplus\cq/\ZZ$ and
$x_it_j=t_jx_i$ for all $1\leq i,j\leq d$.
Then $AC$ becomes an associative algebra.
For $\bff u\in\cd$ set
$$\done u=\sum_{i=1}^d u_i\parder {x_i}+\sum_{i=1}^d u_it_i\parder {t_i}\in\mathrm{Der}(AC).$$
Denote by $\calg$ the Lie subalgebra of $\mathrm{Der}(AC)$ spanned by
$$\{\xone p\done u,\xone l\tonebar s\mid \bff p\in\ndz,\bff l\in\nd,
     \ovbff s\in\Gamma,\bff u\in\cd\}.$$
The Lie algebra $\calg$ has Lie bracket
\begin{equation}\label{eq3.1}
 \begin{aligned}
   &[\xone m\done u,\xone n\done v]=\sum_{i=1}^d u_in_ix^{\bff m+\bff n-\e_i}\done v
                         -\sum_{i=1}^d v_im_ix^{\bff m+\bff n-\e_i}\done u;\\
   &[\xone m\done u,\xone l\tonebar s]=\sum_{i=1}^d u_in_ix^{\bff m+\bff n-\e_i}\tonebar s
                         +(\bff u\mid\bff s)\xtwo ms\tonebar s;\\
   &[\xone p\tonebar r,\xone l\tonebar s]=(\sgmf rs-\sgmf sr)\xtwo ms\tonebar{r+s},
 \end{aligned}
\end{equation}
where $\bff u,\bff v\in\cd,\bff m,\bff n\in\ndz,\bff p,\bff l\in\nd$
and $\bff r,\bff s\in\Gamma_0$.
It is easy to see that $\calg$ is $\Gamma$-graded, i.e.,
$\calg=\bigoplus_{\ovbff s\in\Gamma}\calg(\ovbff s)$, where
$$\begin{aligned}
 &\calg(\ovbff s)=\spanc{\xone n\tonebar s\mid \bff n\in\nd};\\
 &\calg(\ovbff 0)=\spanc{\xone n\tonebar 0,\xone p\done u\mid \bff n\in\nd,
                          \bff p\in\ndz,\bff u\in\cd}.
  \end{aligned}$$
On the other hand, $\calg$ may be decomposed into the sum of two subalgebras
$$\calgx=\spanc{\xone p\done u\mid \bff p\in\ndz,\bff u\in\cd}$$
and
$$\calgt=\spanc{\xone n\tonebar s\mid \bff n\in\nd,\ovbff s\in\Gamma}.$$
Set $\deg{x_i}=1,\ \deg{\parder{x_i}}=-1$ for all $1\leq i\leq d$,
and for $\bff m\in\nd$ write $|\bff m|=\sum_{i=1}^dm_i$.
This gives $\Z$-gradations on the algebras $\calgx$ and $\calgt$.
For $i\in\N$, set
$$\calgx_i=\spanc{\xone p\done u\mid |\bff p|=i+1,\bff u\in\cd},\ \ \ \
  \calgt_i=\spanc{\xone n\tonebar s\mid |\bff n|=i,\ovbff s\in\Gamma}.$$
Then we have $\calgx=\bigoplus_{i\in\N}\calgx_i$
and $\calgt=\bigoplus_{i\in\N}\calgt_i$.

Set $\calg_i=\calgx_i+\calgt_i$, which is a subspace of $\calg$.
These subspaces $\calg_i$ do not make the whole algebra $\calg\ \Z$-graded,
since the derivations $\done u,\ \bff u\in\cd$, are not homogeneous.
However, $\calg_+=\bigoplus_{i\geq 1}\calg_i$
still makes a $\Gamma$-graded ideal of $\calg$.
Write $\calgx_+=\bigoplus_{i\geq 1}\calgx_i$ and $\calgt_+=\bigoplus_{i\geq 1}\calgt_i$.
The main result in this section is the following

\begin{thm}\label{thm3.1}
(1) The commutator $[\calgx,\calgx]=[\calgx_0,\calgx_0]\oplus\calgx_{\geq 1}$;\\
(2) For every finite dimensional $\calgx$-module $U$, there exists an integer $p\gg 0$
    such that $\calgx_pU=0$.\\
(3) For every finite dimensional $\calg$-module $U$, there exists an integer $p\gg 0$
    such that $\calg_pU=0$.\\
(4) The ideal $\calg_+$ annihilates any finite dimensional irreducible $\calg$-module.
\end{thm}
\pf{
(1) is easy to check by equation (\ref{eq3.1}).
For (2) we choose a basis $\mu_1,\cdots,\mu_d$, of $\cd$
with all $\mu_i$ being generic.
Then $\calgx$ has a decomposition $\calgx=\bigoplus_{i=1}^d\calg^{(x,\mu_i)}$,
where $\calg^{(x,\mu_i)}=\spanc{\xone md_{\mu_i}\mid \bff m\in\ndz}$ are subalgebras of $\calgx$.
Notice that for each $1\leq i\leq d$,
the subalgebra $\calg^{(x,\mu_i)}$
is exactly the infinite dimensional Lie algebra $\mathcal L_+$ studied in \cite{BF1}.
By Theorem 3.1 in \cite{BF1} we see that there exists an integer $p_i\gg 0$ such that
$(\calg^{(x,\mu_i)}\cap\calgx_{p_i})U=0$.
Then (2) stands for the integer $p=\max\{p_1,\cdots,p_d\}$.

To prove (3) we consider $U$ as an $\calgx$-module and let $p$ be as in (2).
For $\bff s\in\Gamma_0$ choose $\bff u\in\cd$ such that $(\bff u\mid\bff s)\neq0$.
Then the equation
$$[\xone m\done u,\tonebar s]=(\bff u\mid\bff s)\xone m\tonebar s$$
implies that $(\xone m\tonebar s) U=0$ if $|\bff m|\geq p$.
This proves (3).

Let $V$ be a finite dimensional irreducible $\calg$-module
with representation map $\rho$.
Since $\calg_+$ is an ideal of $\calg$,
the subspace $V'=\{v\in V\mid\calg_+v=0\}$ is a $\calg$-submodule of $V$.
By (3) and equation (\ref{eq3.1})
we see that $\rho(\calg_+)$ is nilpotent.
Then the Lie's Theorem implies that $V'$ is nonzero,
hence $V'=V$ by the irreducibility of $V$.
So $\calg_+V=0$, proving (4).
}

At last we mention that the quotient $\calg/\calg_+$ is isomorphic to
the direct sum $\gld\oplus\gln$ as Lie algebras.
The isomorphism map is given by
\begin{equation}\label{eq3.2}
 \bxone r\mapsto \tonebar r+\calg_+;\ \
    E_{ij}\mapsto  x_id_{\e_j}+\calg_+,
\end{equation}
where $\bff r\in\Gamma_0,1\leq i,j\leq d$ and
$E_{ij}$ is the $d\times d$ matrix with 1 at the $(i,j)$-th entry and 0 elsewhere.
Through this map we also have $\calgx/\calgx_+\cong\gld$ and
$\calgt/\calgt_+\cong\gln$.

%%%%%%%%%%%%%%%%%%%%%%%%%%%%%%%%%%%%%%%%%%%%%%%%%%%%%%%%%%%%%%%%%%%%%%%%%%%%%%% Section 4
\section{Cuspidal $\ZD$-modules}
\def\theequation{4.\arabic{equation}}
\setcounter{equation}{0}

In this section we study cuspidal $\ZD$-modules with support lying in one coset $\al+\zd$
for some $\al\in\cd$.
Fix $M=\bigoplus_{\bff s\in\zd}M_{\al+\bff s}$ such a $\ZD$-module,
where $M_{\al+\bff s}$ is the weight space with weight $\al+\bff s$.
Clearly we have $\tone r M_{\al+\bff s}\subseteq M_{\al+\bff r+\bff s}$.
Since $M$ is a free $\ZZ$-module, we may write
$$M\cong\bigoplus_{\ovbff s}U_{\ovbff s}\otimes\ZZ,$$
where $U_{\ovbff s}=M_{\al+\bff s}$ for $\bff s\in\Gamma_0$.
Set $U=\bigoplus_{\bff s\in\Gamma_0}U_{\ovbff s}$,
which is a $\Gamma$-graded space.

We should introduce some operators in $\U(\D\ltimes\ZZ)$ that act on $U$.
Set $\dtwo um=\tone{-m}\partwo um$ for $\bff u\in\cd$ and $\bff m\inr$,
acts on each $U_{\ovbff s}$, hence on $U$.
For $\bff r\in\Gamma_0, \bff n\inr$, define the operator $\ltwo nr$
to be the restriction of $\tone{-n}\ttwo nr$ 
mapping $U_{\ovbff s}$ to $U_{\bf\ov{r+s}}$ for each $\ovbff s\in\Gamma$.
The commutators among these operators are
\begin{equation}\label{eq4.1}
 \begin{aligned}
   &[\dtwo um,\dtwo vn]=(\bff u\mid\bff n)(D(\bff v,\bff m+\bff n)-\dtwo vn)
                        -(\bff v\mid\bff m)(D(\bff u,\bff m+\bff n)-\dtwo um);\\
   &[\dtwo um,\ltwo ns]=(\bff u\mid\bff n+\bff s)L(\bff m+\bff n,\bff s)
                         -(\bff u\mid\bff n)\ltwo ns;\\
   &[\ltwo mr,\ltwo ns]=(\sgmf rs-\sgmf sr)L(\bff m+\bff n,\bff r+\bff s),
 \end{aligned}
\end{equation}
where $\bff u,\bff v\in\cd, \bff m,\bff n\inr, \bff r,\bff s\in\Gamma_0$.
Since for $v\in U_{\ovbff s},\ \bff n\inr$, we have
\begin{equation}\label{eq4.2}
 \begin{aligned}
  &\partwo um(\tone n v)=(\bff u\mid \bff n)\ttwo mn v+\ttwo mn\dtwo umv;\\
  &\ttwo mr(\tone n v)=\tone n(\ttwo mr v)=\ttwo mn(\ltwo mr v),
 \end{aligned}
\end{equation}
the operators $\dtwo um, \ltwo mr$ completely determine the $\D$-action on $M$.
Furthermore, the operators $\dtwo um$ may be restricted to $U_{\ovbff s}$,
and $\ltwo nr$ may be considered as an operator
from $U_{\ovbff s}$ to $U_{\bf\ov{r+s}}$ for each $\ovbff s\in\Gamma$.
In other words, $\dtwo um$ are of degree $\ovbff 0$
and $\ltwo nr$ are of degree $\ovbff r$ in $\mathrm{End}(U)$.

Let $\mu_1,\cdots,\mu_d$ be a basis of $\cd$ with all $\mu_i$ being generic.
Recall the solenoidal Lie algebra $\g_\mu$ over $\cq$,
the solenoidal Lie algebra $W_\mu$ over the commuting torus $A=\C[x_1^{\pm1},\cdots, x_d^{\pm1}]$.
Then the algebra $\D$ may be decomposed into direct sum of subalgebras
$$\D=\g_{\mu_1}\oplus W_{\mu_2}\oplus\cdots\oplus W_{\mu_d}.$$
Here $W_{\mu_i}$ are solenoidal Lie algebras over the commuting torus $\ZZ$.
By results from \cite{BF1} and \cite{Xu2},
we know that the operators $\dtwo {\mu_i}m, 1\leq i\leq d$ and $\ltwo mr, \bff r\in\Gamma_0$,
have polynomial dependence on $\bff m\inr$ with coefficients in $\mathrm{End}(U)$.
Moreover, the constant term of $\dtwo {\mu_i}m$
on each $U_{\ovbff s}$ is $(\bff u\mid \al+\bff s)Id$.
Hence $\dtwo um$ are polynomials on $\bff m$,
since they are linear combinations of $\dtwo {\mu_i}m$.

Set
$$\dtwo um=\sum_{\bff p\in\nd}\frac{\bff m^{\bff p}}{\bff p!}\ftwo up;\ \ \ \
  \ltwo mr=\sum_{\bff p\in\nd}\frac{\bff m^{\bff p}}{\bff p!}\gtwo rp,$$
where $\ftwo up$ are operators in $\mathrm{End}(U)$ of degree $\ovbff 0$
and $\gtwo rp\in\mathrm{End}(U)$ are operators of degree $\ovbff r$.

Now by expanding the equations in (\ref{eq4.1})
and comparing coefficients at both sides, we get the following commutators
\begin{equation}\label{eq4.3}
 \begin{aligned}
  &[\ftwo up,\ftwo vl]=\begin{cases}
  \sum\limits_{i=1}^du_il_if_{\bff v}({\bf p+l-\e_i})
           -\sum\limits_{i=1}^dv_ip_if_{\bff u}({\bf p+l-\e_i}),
         &\text{ if }\bff p,\bff l\neq \bff 0;\\
  0 ,    &\text{ if }\bff p=\bff 0\text{ or }\bff l=\bff 0,
      \end{cases}\\
  &[\ftwo up,\gtwo sl]=\begin{cases}
         (\bff u\mid\bff s)\gtwo sl, & \text{ if }\bff p=\bff 0;\\
         (\bff u\mid\bff s)g_{\bf s}({\bf p+l})+\sum\limits_{i=1}^du_il_ig_{\bf s}({\bf p+l}-\e_i),
                                      & \text{ if } \bff p\neq \bff 0,
      \end{cases}\\
  &[\gtwo rp,\gtwo sl]=\begin{cases}
  (\sgmf rs-\sgmf sr)g_{\bf{r+s}}({\bf p+l}), & \text{ if }\bff r+\bff s\notinr;\\
          0,                                     & \text{ if }\bff r+\bff s\inr,
      \end{cases}
  \end{aligned}
\end{equation}
where $\bff u,\bff v\in\cd, \bff p,\bff l\in\nd$ and $\bff r,\bff s\in\Gamma_0$.
Recall the Lie bracket of the algebra $\calg$ from equation (\ref{eq3.1})
and we see from equations in (\ref{eq4.3}) that the operators
$$\{\ftwo up,\gtwo sl\mid\bff u\in\cd,\bff s\in\Gamma_0,
       \bff p\in\ndz,\bff l\in\nd\}$$
yield a finite dimensional $\Gamma$-graded representation for $\calg$ on $U$.
Denote by $\rho$ this representation map.
Then combining equation (\ref{eq4.2}) we obtain

\begin{thm}\label{thm4.1}
There exists an equivalence between the category of finite dimensional $\Gamma$-graded
$\calg$-modules and the category of cuspidal $\ZD$-modules with support
lying in some coset $\al+\zd$.
This equivalence functor associates to a finite dimensional $\Gamma$-graded $\calg$-module
$U=\bigoplus_{\ovbff s\in\Gamma}U_{\ovbff s}$ an $\ZD$-module
$$M=\bigoplus_{\bff s\in\zd}U_{\ovbff s}\otimes\tone s$$
with the $\D$-action
\begin{equation}\label{eq4.4}
  \begin{aligned}
    &\partwo um(v_{\ovbff s}\otimes\ttwo sn)=\left((\bff u\mid\al+\bff n+\bff s)\mathrm{Id}+
            \sum_{\bff p\in\ndz}\frac{\bff m^{\bff p}}{\bff p!}\rho(\ftwo up)\right)
            v_{\ovbff s}\otimes\tthree snm;\\
    &\ttwo mr(v_{\ovbff s}\otimes\ttwo sn)=\sum_{\bff p\in\nd}
            \frac{\bff m^{\bff p}}{\bff p!}\rho(\xone p\tonebar r)v_{\ovbff s}\otimes\tfour snmr,
  \end{aligned}
\end{equation}
where $\bff m,\bff n\inr, \bff r\in\Gamma_0\backslash R,
\bff s\in\Gamma_0$ and $v_{\ovbff s}\in U_{\ovbff s}$.
\end{thm}

As a consequence we classify irreducible cuspidal $\ZD$-modules,
which was first done in \cite{LiuZ2}.
First we recall modules of tensor fields over the algebra $\D$.
Let $\al\in\cd$ and let $V$ be a finite dimensional $\gld$-module,
$W=\bigoplus_{\ovbff s\in\Gamma}W_{\ovbff s}$ a finite dimensional $\Gamma$-graded $\gln$-module.
On the tensor space $\bigotimes_{\bff s\in\zd}(V\otimes W_{\ovbff s}\otimes \tone s)$
there is a $\D$-module structure defined by
$$\begin{aligned}
 &\partwo um(v\otimes w\otimes\tone s)=\left((\bff u\mid \al+\bff s)v
                              +(\sum_{i,j=1}^dm_iu_jE_{ij})v\right)\otimes w\otimes\ttwo ms;\\
 &\tone r(v\otimes w\otimes\tone s)=v\otimes (\bxone rw)\otimes\ttwo rs,
\end{aligned}$$
for $\bff u\in\cd, \bff m\inr, \bff r\notinr,\bff s\in\zd, v\in V$ and $w\in W_{\ovbff s}$.
We denote this module by $\vavw$, called a module of tensor fields for $\D$.
Clearly, $\vavw$ becomes a $\ZD$-module provided the $\ZZ$-action
$$\tone n(v\otimes w\otimes\tone s)=v\otimes w\otimes\ttwo ns\text{ for }\bff n\inr.$$
Moreover, $\vavw$ is irreducible as a $\ZD$-module as long as $V, W$ are irreducible.

\begin{thm}\label{thm4.2}
Any irreducible cuspidal $\ZD$-module
is isomorphic to $\vavw$ for some $\al\in\cd$,
finite dimensional irreducible $\gld$-module $V$
and finite dimensional $\Gamma$-graded irreducible $\gln$-module $W$.
\end{thm}

By applying Theorem \ref{thm4.1},
we will prove this theorem again in a rather conceptional way.
Let $M=\bigoplus_{\bff s\in\zd}M_{\al+\bff s}$ be an irreducible cuspidal $\ZD$-module.
By Theorem \ref{thm4.1}, $M\cong M'\otimes \ZZ$,
where $M'=\bigoplus_{\ovbff s\in\Gamma}M_{\al+\bff s}$
is a finite dimensional irreducible $\Gamma$-graded $\calg$-module.
By Theorem \ref{thm3.1}(4) we have $\calg_+M'=0$.
Notice that the ideal $\calg_+$ is still $\Gamma$-graded.
This reduces $M'$ to a $\Gamma$-graded and irreducible module over $\calg/\calg_+$,
which is isomorphic to $\gld\oplus\gln$.
The following lemma is from \cite{Li}.

\begin{lem}\label{lem4.3}
Let $A_1,A_2$ be two associative algebras with identity
and $V$ be an irreducible module over $A_1\otimes A_2$.
Suppose that $A_1$ is of countable dimension.
Then $V$ is isomorphic to an $A_1\otimes A_2$-module of the form $V_1\otimes V_2$,
where $V_1$ is an irreducible $A_1$-submodule of $V$,
and $V_2=\mathrm{Hom}_{A_1}(V_1,V)$ is a natural irreducible $A_2$-module given by
$$(a\cdot f)v=(1\otimes a)(f(v))\text{ for }a\in A_2,f\in\mathrm{Hom}_{A_1}(V_1,V),v\in V_1.$$
\end{lem}

Denote $\gl=\gld\oplus\gln$, which is $\Gamma$-graded with graded spaces
$\gl_{\ovbff 0}=\gld\oplus\C\bxone 0$ and
$\gl_{\ovbff s}=\C\bxone s$ for $\ovbff s\in\Gamma\backslash\{\ovbff 0\}$.
Fix a $\ovbff s\in\Gamma$ such that $M_{\al+\bff s}\neq 0$.
Notice that $M_{\al+\bff s}$ is a $\gld$-module through the $\calgx/\calgx_+$-action.
Fix $V$ an irreducible $\gld$-submodule of $M_{\al+\bff s}$
and set $W_{\ovbff r}=\mathrm{Hom}_{\gld}(V,M_{\al+\bff r+\bff s})$ for any $\ovbff r\in\Gamma$.
Let
$$W=\bigoplus_{\ovbff r\in\Gamma}W_{\ovbff r}
   =\mathrm{Hom}_{\gld}(V,\bigoplus_{\ovbff r\in\Gamma}M_{\al+\bff r+\bff s})
   =\mathrm{Hom}_{\gld}(V,M')$$
 and define
$$(\bxone p\cdot f)v=(1\otimes\bxone p)(f(v))\text{ for }
   v\in V,\ovbff p\in\Gamma, f\in W_{\ovbff r},$$
where $1\otimes\bxone p$ is an element in $\U(\gld)\otimes\U(\gln)=\U(\gld\oplus\gln)$.
Clearly this makes $W$ a $\Gamma$-graded $\gln$-module.
By Lemma \ref{lem4.3} we have $M'\cong V\otimes W$.
Moreover, $W$ is irreducible as a $\gln$-module
since $M'$ is irreducible as a $\gld\oplus\gln$-module.
This proves Theorem \ref{thm4.2}.

From the proof above one can see that
all irreducible $\gld$-submodules of $M'$ are isomorphic to $V$.

\end{document}